\documentclass[11pt]{article}
\usepackage{inputenc}
\usepackage[T2A]{fontenc}
\usepackage[dvips]{graphicx}
\usepackage{amsthm,amsfonts,amsmath,amssymb}
\usepackage{color}
\usepackage{amsmath,bm}

\def\div{\mbox{div}\,}

\begin{document}

\author{\bf A.V. Rukavishnikov }
\title{\bf On the area  of optimal parameters choice for the numerical method of non-stationary hydrodynamics problem with
feature}
\maketitle{}

	\centerline{Computing Center of the Far Eastern Branch of the Russian Academy of Sciences,}
\centerline{680000, Khabarovsk, Russia, e-mail: \tt 78321a@mail.ru}			
\maketitle{}
\begin{abstract}
\begin{tabular}{p{0mm}p{139mm}}
&\noindent {\footnotesize \qquad
For an approximate solution of the non-stationary nonlinear Navier-Stokes equations
for the flow of an incompressible viscous fluid, depending on the set of input data and the geometry of the domain,
the area of optimal parameters in the variables $\nu$ and $\nu^{\ast}$ is experimentally determined depending on $\delta$
included in the definition of the $R_{\nu}$-generalized solution of the problem and the degree of the weight function
 in the basis of the finite element method. To discretize the problem in time, the Runge-Kutta methods
of the first and second orders were used. The areas of optimal parameters for various values
 of the incoming angles are established.

{\bf Keywords:} {nonlinear Navier-Stokes equations, feature, finite element method.}}
\end{tabular}\end{abstract}

\paragraph{Introduction.}

The purpose of this study is to find
the area  of optimal parameters choice for the numerical method of non-stationary Navier-Stokes equations with
feature.
The search for a solution to a nonlinear problem is reduced to a sequence of approximate linear problems using Runge-Kutta methods of 1st and 2nd order. The peculiarity of the study lies in the fact that the domain is a non-convex polygon with a incoming angle at the boundary. Using classical approximate methods, the error arising in the vicinity of the feature point propagates into the inner part of the  domain, where the solution has sufficient smoothness. In this case, the  convergence rate of the approximate solution to the exact one is significantly less than for convex domains. The proposed numerical method overcomes these difficulties. It is based on two ideas, namely, the introduction into the variational formulation of problems a weight function to some order and special basis functions.

In \cite{Ruk2018} the concept of an $R_{\nu}$-generalized solution in weighted sets for the Stokes problem is defined.
The main feature of the variational formulation of this problem, in contrast to the classical formulation  \cite{Brezzi},
is that it is asymmetrical. In \cite{Math2} the existence and uniqueness of an $R_{\nu}$-generalized solution for the Stokes problem is established.

    In the course of numerical experiments, an area  of optimal parameters of the method was determined. The order of convergence of the approximate solution to the exact one of the nonlinear problem is the same for angles more than
    $\pi$ and significantly greater than using classical approaches. Part of the optimal set of free parameters of
    the proposed method doesn't depend on the value of the incoming angle.
    The optimal convergence rate is achieved without using a mesh refinement in the vicinity of the feature point.

\paragraph{1. The problem statement. }

Consider the flow of a viscous incompressible fluid in a 2-dimensional non-convex polygonal
 domain $\Omega$ with an incoming angle $\omega$ on its boundary $\partial \Omega$
with a vertex at the origin ${\cal O}=(0,0)$.
Let ${\bf x}=(x_1, x_2)$ be an element in $R^2,$ $t$ be an element in time and $Q=\Omega\times (0, T)$. Given fields ${\bf u}_0={\bf u}_0({\bf x})$ in $\Omega, {\bf f}={\bf f}({\bf x},t)=\{ f_i({\bf x},t)\}_{i=1}^2$ in $Q$ and ${\bf g}={\bf g}({\bf x},t)=\{g_i({\bf x},t)\}_{i=1}^2$ in $Q$
such that $\int\limits_{\Omega}{\div {\bf g} d{\bf x}}=0$ at each time $t\in(0,T),$ it is required to find the fields
${\bf u}={\bf u}({\bf x},t)=\{u_i({\bf x},t)\}_{i=1}^2$ and $P=P ({\bf x},t)$ such that the following identities hold:
\begin{gather}
 \frac{\partial {\bf u}}{\partial t}-\triangle {\bf u} + \mbox{ curl } {\bf u}\times {\bf u}+ \mbox{ grad } P={\bf f},
 \qquad\, \,\mbox{ div } {\bf w} ={ 0}\, \qquad\mbox{ in } \,\,\qquad Q,
\label{eq:1}\\
\qquad\qquad\,\,\,\,\, {\bf u}({\bf x},0)={\bf u}_0 \qquad \mbox{ in } \qquad\Omega,\qquad\qquad
{\bf u}={\bf g} \qquad \mbox{ on } \qquad \partial \Omega\times (0,T).
\label{eq:2}
\end{gather}

As a time discretization of problem (1)-(2), we use the Runge-Kutta schemes of the 1st and 2nd orders.
To do this, we first introduce the notation ${\bf v}^n={\bf v}^n({\bf x})$ to approximate the function
${\bf v}({\bf x}, n\triangle t), n=0,1,2,...,N$ and ${\bf v}^{n+\gamma}$ to approximate the function
${\bf v}({\bf x}, (n+\gamma)\triangle t), \gamma\in (0,1), n=0,1,2,...,N-1$. Parameter $\triangle t$
is such that $T=N\cdot \triangle t.$
Moreover, let ${\bf v}^{-1}:={\bf v}^0$ and $\bar{\bf v}^{n+1}:=0.5({\bf v}^ {n+1}+{\bf v}^n)$ and ${\bf U}^n$
a suitable approximation to {\bf u} at time $n \triangle t.$

{\bf 1st order scheme}.

Given ${\bf u}^n, {\bf U}^n:=\frac{3}{2} {\bf u}^n -\frac{1}{2} {\bf u}^ {n-1}, P^n, \bar{\bf f}^{n+1}$
  and ${\bf g}^{n+1}$: find ${\bf u}^{n+1}$ and $P^{n+1}$ as a solution to the system of equations:
\begin{gather}
   (\triangle t)^{-1}{\bf u}^{n+1}-\triangle \bar{\bf u}^{n+1} + \mbox{ curl } {\bf U}^n \times \bar{\bf u}^{n+1}+ \mbox{ grad } \bar P^{n+1}=\bar{\bf f}^{n+1}+(\triangle t) ^{-1}{\bf u}^{n}\qquad \mbox{ in } \,\,\qquad \Omega,
   \label{eq:3}\\
\qquad\qquad\qquad\qquad\qquad\qquad\qquad\qquad\qquad\qquad\qquad\qquad \qquad\, \,\mbox{ div } {\bf u}^{n+1} ={ 0} \,\qquad \mbox{ in } \,\,\qquad \Omega,
\label{eq:4}\\
\qquad\qquad\qquad\qquad\qquad\qquad\qquad\qquad\qquad\qquad\qquad\qquad\qquad\qquad{\bf u}^{n+1}={\bf g}^{n+1 } \qquad \mbox{ on } \qquad \partial \Omega.
\label{eq:5}
\end{gather}

{\bf 2nd order scheme}.

This scheme consists of two steps.

{\it Step 1}.

Given ${\bf u}^n, {\bf U}^n:=\frac{3}{2} {\bf u}^n -\frac{1}{2} {\bf u}^ {n-1},
  {\bf f}^{n+\gamma}$ and ${\bf g}^{n+\gamma}$: find ${\bf u}^{n+\gamma}$ and
$P^{n+\gamma} (\gamma\in (0,1))$ as a solution to the system of equations:
\begin{gather}
  (\gamma \triangle t)^{-1}{\bf u}^{n+\gamma}-\triangle {\bf u}^{n+\gamma} + \mbox{ curl } {\bf U}^n\times {\bf u}^{n+\gamma}+ \mbox{ grad } P^{n+\gamma}={\bf f}^{n+\gamma}+(\gamma\triangle t)^{-1}{\bf u}^{n}\qquad \mbox{ in } \,\,\qquad \Omega,
  \label{eq:6}\\
 \qquad\qquad\qquad\qquad\qquad\qquad\qquad\qquad\qquad\qquad\qquad\qquad\qquad\, \,\mbox{ div } {\bf u}^{n+\gamma} ={ 0}\,\qquad \mbox{ in } \,\,\qquad \Omega,
\label{eq:7}\\
\qquad\qquad\qquad\qquad\qquad\qquad\qquad\qquad\qquad\qquad\qquad\qquad\qquad\qquad{\bf u}^{n+\gamma}={\bf g}^{n+\gamma} \qquad \mbox{ on } \qquad \partial \Omega.
\label{eq:8}
\end{gather}

{\it Step 2}.

Given ${\bf u}^n, {\bf u}^{n+\gamma}, P^{n+\gamma},
{\bf U}^n:=\frac{3}{2} {\bf u}^n -\frac{1}{2} {\bf u}^{n-1}, {\bf f} ^{n+1},
{\bf f}^{n+\gamma}$ and ${\bf g}^{n+1}$: find ${\bf u}^{n+1}$ and $P^{n+1}$ as a solution to the system of equations:
\begin{gather}
  (\triangle t)^{-1}{\bf u}^{n+1}+\gamma (-\triangle {\bf u}^{n+1} + \mbox{ curl } {\bf U}^n\times {\bf u}^{n+1}+ \mbox{ grad }  P^{n+1})=\\\nonumber
    =(\triangle t)^{-1}{\bf u}^{n}+\gamma {\bf f}^{n+1}+(1-\gamma){\bf f}^{n+\gamma}
  -(1-\gamma) (-\triangle {\bf u}^{n+\gamma} + \mbox{ curl } {\bf U}^n\times {\bf u}^{n+\gamma}+ \mbox{ grad } P^{n+\gamma})
  \quad \mbox{ in } \quad \Omega,
  \label{eq:9}\\
\qquad\qquad\qquad\qquad\qquad\qquad\qquad\qquad\qquad\qquad\qquad\qquad\qquad \qquad\, \,\mbox{ div } {\bf u}^{n+1} ={ 0}\,\qquad \mbox{ in } \,\,\qquad \Omega,
\label{eq:10}\\
\qquad\qquad\qquad\quad\qquad\qquad\qquad\qquad\qquad\qquad\qquad\qquad\qquad\qquad\qquad{\bf u}^{n+1}={\bf g}^{n+1} \qquad \mbox{ on } \qquad \partial \Omega.
\label{eq:11}
\end{gather}

At each step of both schemes, it is necessary to be able to solve the following problem: find the fields ${\bf v}=(v_1,v_2)$ and $q$ such that
\begin{gather}
\theta {\bf v}-\triangle {\bf v} + W \times {\bf v}+ \mbox{ grad } q={\bf F} \qquad\mbox{ in } \,\,\qquad \Omega,
\label{eq:12}\\
 \qquad\qquad\qquad\qquad\qquad\, \,\mbox{ div } {\bf v} ={ 0}\, \qquad\mbox{ in } \,\,\qquad \Omega,
\label{eq:13}\\
\qquad\qquad\qquad\qquad\qquad\qquad
{\bf v}={\bf G} \qquad \mbox{ on } \qquad \partial \Omega,
\label{eq:14}
\end{gather}
where ${\bf F}$ and $W$ are given functions on $\Omega$ and ${\bf G}$ is given on $\partial \Omega.$

Let us define an $R_{\nu}$-generalized solution of problem (12)-(14) in the domain $\Omega,$
  having an incoming angle on the $\partial \Omega$ with a vertex at the origin.
To do this, we define the necessary weight sets. First, we introduce the concept of a weight function
$\rho({\bf x}): \rho({\bf x})=\{\sqrt{x_1^2+x^2_2}, \mbox{ if }  {\bf x}\in \Omega_{\delta}
\mbox{ and } \delta, \mbox{ if } {\bf x}\in \bar{\Omega}\diagdown \Omega_{\delta}\},$ where
 $\Omega_{\delta}=\{{\bf x}\in \bar{\Omega}: \sqrt{x_1^2+x^2_2}\leq \delta\}$ and $\delta\ll 1.$

Denote by $L_{2,\alpha}(\Omega,\delta)$ the set of functions $s({\bf x})$ satisfying the conditions:\\
1) $\|\rho^{\alpha}s\|_{L_2(\Omega\diagdown \Omega_{\delta})}\geq C_1>0;$\\
2) $|s({\bf x})|\leq C_2 \delta^{\alpha-\varepsilon} \rho^{\varepsilon-\alpha}({\bf x}),{\bf x}\in \Omega_{\delta},$\\
where $C_2$ is a positive constant independents $s({\bf x})$, $\varepsilon$ is a small positive parameter independent of
$\delta, \alpha, s({\bf x}),$ with bounded norm
$\|s\|_{L_{2,\alpha}(\Omega)}:=\|\rho^{\alpha} s\|_{L_2(\Omega)}$ of the space
$L_{2,\alpha}(\Omega).$ $L^0_{2,\alpha}(\Omega,\delta)$ subset of $L_{2,\alpha}(\Omega)$
  such that $s\in L^0_{2,\alpha}(\Omega,\delta),$ if $s\in L_{2,\alpha}(\Omega)$ and $\|\rho^{ \alpha}s\|_{L_1(\Omega)}=0.$

Denote by $W^1_{2,\alpha}(\Omega,\delta)$ the set of functions $s({\bf x})$ satisfying the conditions:\\
1) $\|\rho^{\alpha}s\|_{L_2(\Omega\diagdown \Omega_{\delta})}\geq C_1>0;$\\
2) $|s({\bf x})|\leq C_2 \delta^{\alpha-\varepsilon} \rho^{\varepsilon-\alpha}({\bf x}),{\bf x}\in \Omega_{\delta},$\\
3) $|D^1 s({\bf x})|\leq C_2 \delta^{\alpha-\varepsilon} \rho^{\varepsilon-\alpha-1}({\bf x}), { \bf x}\in \Omega_{\delta},$\\
where $C_2$ is a positive constant independent $s({\bf x})$, $\varepsilon$ is a small positive parameter,
independent of $\delta, \alpha, s({\bf x}),$ with bounded norm
  $\|s\|_{W^1_{2,\alpha}(\Omega)}:=\sqrt{\sum\limits_{|k|\leq 1}{\|\rho^{\alpha} | D^k s|\|^2_{L_2(\Omega)}}}$ of the space
$W^1_{2,\alpha}(\Omega).$ $\hat W^1_{2,\alpha}(\Omega,\delta)$ subset
$W^1_{2,\alpha}(\Omega)$ such that $s\in \hat W^1_{2,\alpha}(\Omega,\delta),$
if $s\in W^1_{2,\alpha}(\Omega, \delta)$ and $ s=0$ on $\partial \Omega.$

Let $L_{\infty,\alpha}(\Omega,C_3)$ be the set of functions $s({\bf x})$ with a norm
$\|s\|_{L_{\infty,\alpha}(\Omega,C_3)}= \mbox{vrai} \max\limits_{{\bf x}\in \Omega}
 {|\rho^{\alpha}({\bf x})s({\bf x})|}\leq C_3,
$
where $C_3>0$ independent of $s({\bf x})$.

Let us define an $R_{\nu}$-generalized solution of problem (12)-(14).

{\bf Definition 1}.
{\it Pair $({\bf v}_{\nu}, q_{\nu})\in {\bf W}^{1}_{2,\nu}(\Omega, \delta)\times L^0_{2,\nu}(\Omega,\delta)$ is called
$R_{\nu}$-generalized solution of problem (12)-(14) if ${\bf w}_{\nu}$ satisfies (14) on $\partial \Omega$ and
identities }
 \begin{gather}
a({\bf v}_{\nu},{\bf z})+b({\bf z}, q_{\nu})=l({\bf z}),
\label{eq:15}\\
 c({\bf v}_{\nu},s)=0\qquad\qquad\qquad
\label{eq:16}
\end{gather}
{\it hold, for all $({\bf z}, s)\in \hat{\bf W}^{1}_{2,\nu}(\Omega, \delta)\times L^0_{2,\nu}(\Omega,\delta)$. Here
$({\bf F},W,{\bf G})\in {\bf L}_{2,\alpha}(\Omega,\delta)\times L_{\infty,\beta}(\Omega,C_3)\times {\bf W}^{1/2}_{2,\alpha}(\partial \Omega, \delta), \nu\geq\alpha\geq 0, \beta\leq 2.$ }

We have
$
a({\bf w},{\bf z})=\int\limits_{\Omega}{\Bigl[\theta {\bf w}\cdot (\rho^{2 \nu} {\bf z})+\nabla{\bf w}:\nabla(\rho^{2 \nu} {\bf z})+(W\times {\bf w})\cdot (\rho^{2 \nu} {\bf z})\Bigr] d{\bf x}},
$\\
$
b({\bf z}, p)=-\int\limits_{\Omega}{p \mbox { div } (\rho^{2 \nu} {\bf z})}d{\bf x},\,\,\,
c({\bf w},s)=-\int\limits_{\Omega}{(\rho^{2 \nu} s) \mbox{ div } {\bf w}} d{\bf x},\,\,\,
l({\bf z})=\int\limits_{\Omega}{{\bf F} \cdot (\rho^{2\nu} {\bf z})} d{\bf x}.
$

 {\bf Remark 1}. {\it $b(\cdot,\cdot)\neq c(\cdot,\cdot)$.
Consequently, the variational problem (15), (16) is not symmetric, in contrast to the standard setting (see \cite{Brezzi}).}

{\bf Remark 2}. {\it If $W\in L_{\infty,\beta}(\Omega, C_3), \beta\leq 2, {\bf g=0} $ on $\partial \Omega,$
then there exists a unique $R_{\nu}$-generalized solution $({\bf w}_{\nu},q_{\nu})$
problem (12)-(14) in the asymmetric formulation (15), (16) (see Theorem 5 in \cite{RukJCAM2023}).}

{\bf Remark 3}. {\it Scheme 2 (6)-(11) can be applied when $\gamma $
  equal to $1-\frac{\sqrt{2}}{2}$ (strongly L-stable method \cite{Alex}).
Scheme 1 can be applied due to the validity of Theorem 6 \cite{RukJCAM2023}.}

\paragraph{2. Creation of an approximate approach.}

We will build a quasi-uniform fragmentation $I^h$ of the domain $\bar{\Omega}$ into triangles $K_i$. Their sides of order $h$ which are the essential elements. We split each of them into 3 using the center of mass, which are finite elements $L_{i_j}$ that make up the fragmentation $J^h$.

We  define the main finite element spaces.

1. {\it For the velocity components}. We will use Lagrangian elements of the 2nd order with nodes at the vertices and midpoints of the sides of $L_{i_j}$. The linear span of basic functions $\chi_k({\bf x})$  will be denote by $W_h$.

2. {\it For the pressure}. As approximation nodes, we use the vertices of finite elements, and the vertices of neighboring finite elements are different nodes. On each $L_{i_j}$ we define 1st order basis functions $\theta_l({\bf x})$ whose
support is only one finite element. The linear span of such basis functions forms  a space $S_h$, consisting of functions discontinuous in the domain under consideration.

The pair ${\bf W}_h \times S_h ({\bf W}_h=W_h \times W_h)$ is a Scott-Vogelius 2nd order one \cite{Scott}.

Let's multiply the basis functions of the spaces $W_h$ and $S_h$ by the weight function $\rho({\bf x})$ in some powers ($-\nu^{\ast}$) and (-$\mu^{\ast}$). The values of the powers will be determined later. We define new basis functions
$\phi_k({\bf x})=\chi_k({\bf x}) \cdot \rho^{-\nu^{\ast}}({\bf x}), \psi_l({\bf x})=\theta_l({\bf x}) \cdot \rho^{-\mu^{\ast}}({\bf x})$.

Their linear spans form finite-dimensional spaces $V_h$ and $Q_h$ respectively. $\hat V_h=\{u_h\in V_h: u_h(M_i)=0 \mbox{ where } M_i \mbox{ are nodes on } \partial \Omega\}$. Having found the solution $\hat v_i$ and
$\hat q_j$ at the nodes $M_i$ and $N_j$
for the velocity components and pressure of system of equations (presented below) it is necessary to restore
the true values at the nodes $M_i$ and $N_j$ using the formulas $v_i=\hat v_i \cdot\rho^{-\nu^{\ast}}(M_i)$ and $q_j=\hat q_j \cdot\rho^{-\mu^{\ast}}(N_j).$
We have ${\bf V}_h=V_h\times V_h\subset {\bf W}^1_{2,\nu}(\Omega, \delta), \hat{\bf V}_h=\hat V_h\times \hat V_h\subset \hat{\bf W}^{1}_{2,\nu}(\Omega, \delta), Q_h\subset L^0_{2,\nu}(\Omega, \delta).$

We are all set to determine an approximate $R_{\nu}$-generalized solution of the problem (12)-(14).

{\bf Definition 2}. {\it We will say that a pair of functions
$ ({\bf v}^h_{\nu}, \, q^h_{\nu})$ from the spaces $ {\bf V}_{h} \times Q_{h}, $
satisfying condition (14) at the nodes on $\partial \Omega$, is
an approximate $ R_{\nu} $-generalized solution of the problem (12)-(14)
if for all  pairs of functions
  $ ({\bf z}^h, s^h)$ from the spaces $\hat{\bf V}_{h} \times Q_{h} $
  the following relations
\begin{gather}
a({\bf v}^h_{\nu}, {\bf z}^h)+b({\bf z}^h, q^h_{\nu})=l({\bf z}^h),
\label{eq:17}\\
 c({\bf v}^h_{\nu}, s^h)=0\qquad\quad\,\,\,\qquad\quad
\label{eq:18}
\end{gather}
 hold}.

{\bf Remark 4}.

{\it How to solve the system (17), (18) by the iterative method see \cite{Ruk_Vich Tech}.}

{\bf Remark 5}.

{\it If $ \nu = \nu^{\ast} = \mu^{\ast} = 0 $ then we have an  approximate generalized solution $ ({\bf v}^h, \, q^h) $
of the problem (12)-(14).}

\begin{figure}[h!]
\begin{center}
\includegraphics[scale=0.55]{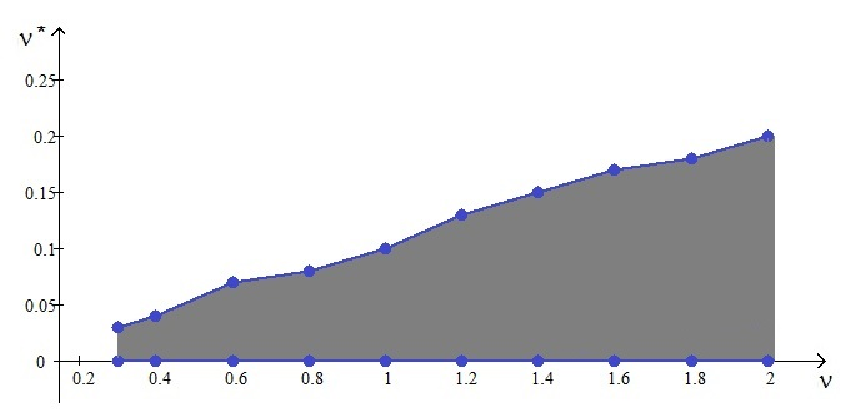}
\caption{Optimal parameters of the weighted FEM for $(\nu, \nu^{\ast}), \delta\in[0.025, 0.035], \omega=\frac{9 \pi}{8}$.}
\end{center}
\end{figure}

\paragraph{3. The results of numerical experiments.}

Let's carry out a number of numerical experiments to find an approximate solution to problem (1)-(2)
 as a sequence of solving problem (12)-(14) of both formulation schemes (17)-(18).
  Consider domains $\Omega_k, k=1,2,3$ with incoming angle $\omega_k,$ where
$$
\Omega_0=\{(x_1, x_2): -1<x_1<1, 0<x_2<1\},
$$
$$
\bar{\Omega}_1=\bar{\Omega}_0\cup \{(x_1, x_2): -1 \leq x_1\leq 0, -1\leq x_2\leq 0\},
$$
$$
\bar{\Omega}_2=\bar{\Omega}_0\cup \{(x_1, x_2): -1 \leq x_1\leq 0, x_1\leq x_2\leq 0\},
$$
$$
\bar{\Omega}_3=\bar{\Omega}_0\cup \{(x_1, x_2): -1 \leq x_1\leq 0, \frac{1}{2} x_1\leq x_2\leq 0\}.
$$
In such cases $\omega_k=(1+2^{-k}) \pi, k=1,2,3.$

Denote by ${\bf u}^{h_i}_{\nu}$ and ${\bf u}^{h_i}$ an approximate $R_{\nu}$-generalized and generalized solutions
(velocity field) of the solution of the problem (1)-(2) at each moment of time.
In the second case $(\nu=\nu^{\ast}=\mu^{\ast}=0, \delta=1).$
In the first case $(\nu,\nu^{\ast},\mu^{\ast}, \delta)$ is the set of free parameters of the weighted finite element method.

The exact solution of the problem (1)-(2) of the velocity and pressure fields depend
on the value of the incoming angle $\omega_k$ in each
moment of time do not belong to the spaces ${\bf W}^2_2(\Omega_k)$ and ${\bf W}^1_2(\Omega_k)$, respectively, and
have the form in polar coordinates $(r,\theta):$\\
$
u_1(r,\theta,t)=e^t\Bigl(r^{\lambda_k}\chi_1(\theta)+\psi_1(r,\theta) \Bigr),
$
$
u_2(r,\theta,t)=e^t\Bigl(r^{\lambda_k}\chi_2(\theta)+\psi_2(r,\theta) \Bigr),
$
$
P(r,\theta,t)=e^t  r^{\lambda_k-1} \gamma(\theta),
$\\
where $\psi_i(r, \theta)$ is the regular part of $u_i(r,\theta,t),$
i.e. a function belonging to the space ${\bf W}^2_2(\Omega_k),$ and
$r^{\lambda_k}\chi_i(\theta)$ and $r^{\lambda_k-1} \gamma(\theta)$ are the singular solution parts
 of the velocity and pressure fields. The exponent $\lambda_k$
is such that it coincides with the smallest real positive solution of the equation
 $\sin(\lambda \omega_k)+\lambda \sin \omega_k=0$.
Thus $(\lambda_1,\lambda_2,\lambda_3)$ take the following approximate values $(0.5445, 0.6736, 0.8008).$

We have
$$
\chi_1(\theta)=\cos(\theta)\, \Xi'_k(\theta)+(1+\lambda_k)\, \sin(\theta)\,  \Xi_k(\theta),
$$
$$
\chi_2(\theta)=(\lambda_k-1)\, \cos(\theta)\, \Xi_k(\theta)+ \sin(\theta)\,  \Xi'_k(\theta),
$$
$$
\gamma(\theta)=(\lambda_k-1)^{-1}\, \Bigl(  \Xi'''_k(\theta)+ (1+\lambda_k)^2\,  \Xi'_k(\theta)\Bigr).
$$

Here
$$
\Xi_k(\theta)=\Bigl[(1+\lambda_k)^{-1} \sin((1+\lambda_k)\theta)-(1-\lambda_k)^{-1} \sin((1-\lambda_k)\theta) \Bigr]
\cos(\lambda_k \omega_k)+\cos((1-\lambda_k)\theta)-\cos((1+\lambda_k)\theta),
$$
$\Xi'_k(\theta)$ and $\Xi'''_k(\theta)$ are its the first and third derivatives with respect to $\theta$, respectively.

\begin{figure}[h!]
\begin{center}
\includegraphics[scale=0.55]{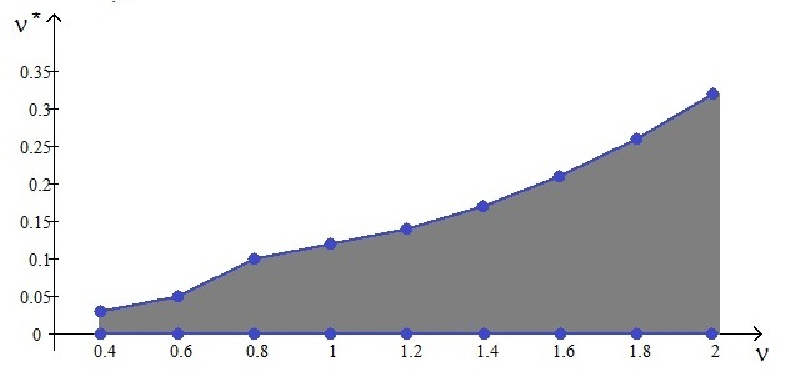}
\caption{Optimal parameters of the weighted FEM for $(\nu, \nu^{\ast}), \delta\in[0.025, 0.035], \omega=\frac{5 \pi}{4}$.}
\end{center}
\end{figure}

In test cases, consider different steps $h_j=2^{1-j}\cdot h, h=0.025, j=1,2,3.$
 Time step $\triangle t=0.01, T=0.5.$ Shown earlier as for stationary \cite{Vich_Mex22}
 and non-stationary problems \cite{RukJCAM2023}
problems, it is possible to determine sets of optimal parameters for which the order of convergence
of the approximate solution to
exact solution is equal to ${\cal O}(h),$ independent of the incoming angle $\omega_k$
in the norm ${\bf W}^1_{2,\nu}(\Omega_k).$ To determine the range of choice of optimal approach parameters,
 we fix the range $\delta \sim h_1: \delta\in[0.025, 0.035].$ Let $\nu^{\ast}=\mu^{\ast}$ and will take non-negative values.
Moreover, each $\omega_k$ will have its own range of value change. The parameter $\nu$ is positive and not exceed
the value  2.

For the first scheme, consider the case
$$
\psi_1(r,\theta)=\psi_2(r,\theta)=0
$$
(the solution contains only singular components).

For the second scheme
$$\psi_1(r,\theta)=\sin (r \cos(\theta))\cdot \cos(r \cos(\theta)),
\psi_2(r,\theta)=- \cos (r \cos (\theta))\cdot \sin (r \sin (\theta)).$$
We will assume that the set $(\nu^{\ast}, \nu)$ falls
into the area of choice of optimal parameters of the numerical method for solving problem (1)-(2),
if it differs by no more than 5 percent from the optimal value in terms of convergence in each
moment of time for all $h_j, j=1,2,3.$
The Figures  1-3 show the area of choice of optimal parameters in the corresponding ranges
for the first scheme for all values of the angle $\omega_k, k=1,2,3.$ For the second scheme, the results similar in structure.

\begin{figure}[h]
\begin{center}
\includegraphics[scale=0.55]{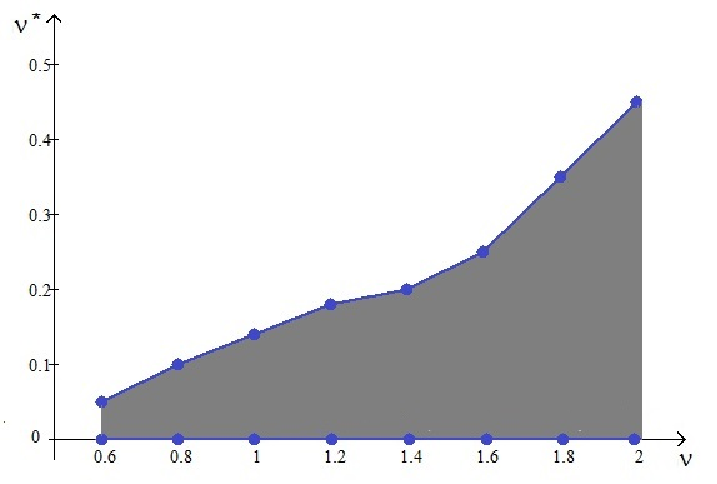}
\caption{Optimal parameters of the weighted FEM for $(\nu, \nu^{\ast}), \delta\in[0.025, 0.035], \omega=\frac{3 \pi}{2}$.}
\end{center}
\end{figure}

\paragraph{Conclusions.}

It is also necessary to investigate other values of the incoming angle $\omega$ and determine the intersection
areas of choice the optimal parameters of the method for each $\omega$. In order to
 establish the area of optimal parameters of the method
for which the required order of convergence is guaranteed regardless   of the incoming angle value.

\paragraph{Acknowledgments.}

The reported study was supported by Russian Science Foundation, project No. 21-11-00039,
https://rscf.ru/en/project/21-11-00039/. The results were obtained using the equipment of SRC ”Far Eastern
Computing Resource” IACP FEB RAS (https://cc.dvo.ru).

\end{document}